\documentclass[fleqn]{mat01}
\usepackage{times,mathtimy,amssymb,latexsym}
\begin{document}

\setcounter{page}{445}
\firstpage{445}

\def\d{\mbox{\rm d}}

\newtheorem{theore}{Theorem}
\renewcommand\thetheore{\arabic{section}.\arabic{theore}}
\newtheorem{theor}{\bf Theorem}
\newtheorem{rem}[theore]{Remark}
\newtheorem{propo}[theore]{\rm PROPOSITION}
\newtheorem{lem}[theore]{Lemma}
\newtheorem{definit}[theore]{\rm DEFINITION}
\newtheorem{coro}[theore]{\rm COROLLARY}
\newtheorem{exampl}[theore]{Example}
\newtheorem{case}{Case}

\renewcommand\thecase{({\it \roman{case}})}

\def\corol{\trivlist \item[\hskip \labelsep{COROLLARY.}]}
\def\noteproof{\trivlist \item[\hskip \labelsep{\it Note added in Proof.}]}

\renewcommand{\theequation}{\thesection\arabic{equation}}

\title{Vector bundles with a fixed determinant on an irreducible nodal curve}

\markboth{Usha N Bhosle}{Vector bundles with a fixed determinant}

\author{USHA N BHOSLE}

\address{Tata Institute of Fundamental Research,
Homi Bhabha Road, Mumbai~400~005, India\\
\noindent E-mail: usha@math.tifr.res.in}

\volume{115}

\mon{November}

\parts{4}

\pubyear{2005}

\Date{MS received 7 July 2005; revised 31 August 2005}

\begin{abstract}
Let $M$ be the moduli space of generalized parabolic bundles
(GPBs) of rank $r$ and degree $d$ on a smooth curve $X$. Let
$M_{\bar L}$ be the closure of its  subset consisting of GPBs
with fixed determinant ${\bar L}$. We define a moduli functor
for which $M_{\bar L}$ is the coarse moduli scheme. Using the
correspondence between GPBs on $X$ and torsion-free sheaves on a
nodal curve $Y$ of which $X$ is a desingularization, we show that
$M_{\bar L}$ can be regarded as the compactified moduli
scheme of vector bundles on $Y$ with fixed determinant. We get a
natural scheme structure on the closure of the subset consisting
of torsion-free sheaves with a fixed determinant in the moduli
space of torsion-free sheaves on $Y$. The relation to
Seshadri--Nagaraj conjecture is studied.
\end{abstract}

\keyword{Nodal curves; torsion-free sheaves; fixed determinant.}

\maketitle

\section{Introduction}

Generalized parabolic vector bundles (GPBs) on a smooth curve $X$
are vector bundles on $X$ together with parabolic structures on
finitely many disjoint divisors $D_j, \, j= 1,\dots, m$
\cite{1,2}. There is an open subscheme $M''$ of the moduli space
$M$ of GPBs on which one can define a determinant morphism into
the moduli space of generalized parabolic line bundles ${\bar
L}$, the map does not extend to $M$.  Let $M''_{\bar L}$ be
its locally closed subset consisting of GPBs with a fixed
determinant ${\bar L}$. In this note, we define a moduli
functor and construct a coarse moduli scheme $M_{\bar L}$ for
it. The moduli scheme contains  $M''_{\bar L}$ as an open dense subscheme.

Let $Y$ be an irreducible projective nodal curve with nodes $y_j,
j= 1, \dots, m$ and $p\hbox{\rm :}\ X \rightarrow Y$ its
desingularization with $D_j$ the inverse image of $y_j$. Denote by
$U$ the moduli variety of torsion-free sheaves of rank $r$, degree
$d$ on $Y$. Let $U'$ be the open subvariety of $U$ corresponding
to vector bundles on $Y$. There is a surjective morphism $f$ from
$M$ onto $U$ \cite{1,2}. The restriction of the morphism $f$ to
$M'= f^{-1}U'$ is an isomorphism onto the open subvariety $U'$ of
$U$. A GPB ${\bar L}$ gives a torsion-free sheaf ${\cal L}$
on $Y$. If ${\cal L}$ is locally free, let $U'_{\cal L}$ be the
closed subset of $U'$ consisting of vector bundles with fixed
determinant ${\cal L}$. Using $f, \, U'_{\cal L}$ may be
identified with $M'_{\bar L}$ and $M_{\bar L}$ can be
regarded as compactified moduli variety of vector bundles on $Y$
with determinant ${\cal L}$.

We show that $f(M_{\bar L}) =\overline{{U'_{\cal L}}}$, the
closure of $U'_{\cal L}$ in $U$, thus giving $\overline{{U'_{\cal
L}}}$ the scheme structure of an image subscheme of $U$. Let $I_j$
denote the ideal sheaf at the node $y_j$. For a torsion-free sheaf
$F$ of rank $r$ on $Y$, let $N= \Lambda^r F/\hbox{(torsion)}$
where (torsion) denotes the torsion subsheaf. Then we show that
for any ${\bar L}$, the image $f(M_{\bar L})$ can be
described (as a set) by
\begin{equation*}
f(M_{\bar L}) = \{F\in U\hbox{\rm :}\  I_j^r{\cal L} \subset
N \subset {\cal L}, \quad \forall j \}.
\end{equation*}
This gives a proof of a conjecture by Seshadri and Nagaraj
(Conjecture (a), p.~136 of \cite{3}). Proving Seshadri--Nagaraj
conjecture was not the aim of this note. The conjecture was proved
by Sun \cite{6} by degeneration methods. However he does not get a
scheme structure on $\overline{{U'_{\cal L}}}$ or a moduli
functor (except in some low rank cases). Our aim is to give a
moduli functor and an explicit construction of a projective moduli
space for it which contains an open subvariety isomorphic to
$U'_{\cal L}$ if ${\cal L}$ is a line bundle. We also deal with
the case when ${\cal L}$ is torsion-free but not locally free. The
construction is much simpler than that of Schmitt \cite{4} and
hence the moduli space is easier to study. For example, properties
like reduced, irreducible, Cohen--Macaulay follow immediately for
our moduli space. Normality is true in rank $2$ and is expected to
be true in general. These properties have been used in computation
of Picard groups in the rank two case.

\section{The moduli scheme of GPBs with fixed determinant}

\renewcommand\thesubsection{\thesection\arabic{subsection}.}

\subsection{}

Let $X$ be a nonsingular projective curve over an algebraically
closed base field $k$. Let $D_j, \, j=1, \dots, m$ be disjoint
divisors on $X$ with $D_j = x_j+x'_j$, where $x_j, x'_j$ are
distinct closed points. We recall here some basics on generalized
parabolic bundles (GPBs), details may be found in \cite{1,2}.

\begin{definit}$\left.\right.$\vspace{.5pc}

{\rm  \noindent A generalized parabolic  bundle (GPB, in short) of
rank $r$ and degree $d$ on  $X$ is a vector bundle $E$ of rank $r$
and degree $d$ on $X$ together with  $r$-dimensional vector
subspaces $F_j(E)$ of $E_{x_j} \oplus E_{x'_j}$. For a subbundle
$N$ of $E$, define $F_j(N)= F_j(E) \cap (N_{x_j}\oplus N_{x'_j})$
and $f_j(N)= \dim F_j(N)$.}
\end{definit}

\begin{definit}$\left.\right.$\vspace{.5pc}

{\rm  \noindent Fix a rational number $\alpha \in (0,1]$. A GPB
$(E,F_j(E))$ is $\alpha$-stable (resp. $\alpha$-semistable) if for
every proper subbundle $N$ of $E$, one has $(d(N)+\alpha \Sigma _j
f_j(N))/r(N)< ({\rm resp.} \leq)$ $(d(E)+\alpha r m)/r.$}
\end{definit}

\begin{definit}$\left.\right.$\vspace{.5pc}

{\rm  \noindent Let $p_j\hbox{\rm :}\ F_j(E) \rightarrow E_{x_j},
\, {p'}_{\!\!j}\hbox{\rm :}\ F_j(E) \rightarrow E_{x'_j}$ be the
projections. Assume that for each $j$, at least one of $p_j, p'_j$
is an isomorphism. The subspace $F_j(E)$ determines an element
$F_j(E)$ of $\hbox{Gr}(r, E_{x_j}\oplus E_{x'_j}) \subset {\bf
P}(\Lambda^r (E_{x_j}\oplus E_{x'_j}))$. One has a (rational)
morphism $\delta\hbox{\rm :}\ {\bf P}(\Lambda^r (E_{x_j}\oplus
E_{x'_j})) \rightarrow {\bf P}(\Lambda^r E_{x_j}\oplus \Lambda^r
E_{x'_j})$. Let $\det F_j(E)$ denote the one-dimensional subspace
of $\Lambda^r E_{x_j}\oplus \Lambda^r E_{x'_j}$ determined by
$\delta(F_j(E))$. Define the determinant of $(E, F_j(E))$ to be
the generalized parabolic line bundle $(\det E, \det F_j(E))$.}
\end{definit}

\begin{definit}$\left.\right.$\vspace{.5pc}

{\rm  \noindent A family of GPBs of rank $r$, degree $d$
parametrized by a scheme $T$ is a tuple $({\cal E}, F_j({\cal
E})_j)$ where ${\cal E} \rightarrow T\times X$ is a family of
vector bundles of rank $r$, degree $d$ on $X$ which is flat over
$T$ and $F_j({\cal E})$ is a rank $r$ subbundle of ${\cal
E}\mid_{T\times {x_j}}\oplus \,  {\cal E}\mid_{T\times {x'_j}}$.
The notion of equivalence of families is the obvious one.

We fix a generalized parabolic line bundle ${\bar L} := (L,
F_j(L))$. Fix isomorphisms $h_j\hbox{\rm :}\ L_{x_j} \rightarrow
k, \, h'_j\hbox{\rm :}\ L_{x'_j} \rightarrow k$. Then $F_j(L)$ can
be identified to a point $F_j(L)$ of ${\bf P}^1$ of the form
$(1:0), (0:1)$ or $(1 \ \hbox{\rm :}\ \lambda_j), \, \lambda_j \in
k^*$.}
\end{definit}

\renewcommand\thesubsection{\thesection\arabic{subsection}}
\subsection{\it The moduli functor}

For simplicity, let us assume that there is only one divisor $D=
x_1+x_2$. Let $({\cal E}, F({\cal E})) \rightarrow T\times X$ be a
family of GPBs of rank $r$, degree $d$ on $X$ with ${\cal E}_t,
t\in T,$ of fixed determinant $L$. For $i=1,2$ we have vector
bundles
\begin{equation*}
{\cal E}_{x_i}= {\cal E}\!\mid_{T\times {x_i}} \rightarrow T.
\end{equation*}
Let $ {\cal G}r \rightarrow T$ denote the Grassmannian bundle of
rank $r$ subbundles of ${\cal E}_{x_1}\oplus {\cal E}_{x_2}$. It
is embedded as a closed subvariety in ${\bf P} (\Lambda^r ({\cal
E}_{x_1}\oplus {\cal E}_{x_2}))$ by Pl\"ucker embedding. Note that
$F({\cal E})$ defines a section of ${\cal G}r$. Since $\det {\cal
E}\!\mid_{t\times X} = L$, it follows that $\det {\cal E} =
{p_T}^*N\otimes p_X^*L$ for some line bundle $N$ on $T$. Hence for
$i= 1, 2$ one has $\det {\cal E}_{x_i}= N \otimes L_{x_i} = N$,
using the isomorphism $h_i$. Let $q_i\hbox{\rm :}\ \Lambda^r
({\cal E}_{x_1}\oplus{\cal E}_{x_2}) \rightarrow \det {\cal
E}_{x_i} = N$ be the projections, $i=1, 2$. Define a hyperplane
subbundle ${\cal H}$ of ${\bf P} (\Lambda^r({\cal
E}_{x_1}\oplus{\cal E}_{x_2}))$ by $q_2= 0$ if $F(L)=(1:0)$, by
$q_1=0$ if $F(L)=(0:1)$ and by $q_2 - \lambda_j q_1=0$ if
$F(L)=(1:\lambda)$. Let $H_T  := {\cal G}r \cap {\cal
H}$. It is a closed reduced subscheme of ${\cal G}r \subset {\bf
P} (\Lambda^r({\cal E}_{x_1}\oplus {\cal E}_{x_2}$)). Note that
$H_T$ is independent of the choice of $h_1, h_2$.

More generally, if we consider parabolic structures over finitely
many disjoint divisors $D_j = x_j + x'_j$, for each $j$ one
constructs the hyperplane bundle ${\cal H}_j$ and Grassmannian
bundle ${\cal G}r_j$ over $T$. Let ${\cal G}r$ be the fibre
product of ${\cal G}r_j$ over $T$ and $H_T$ the fibre product of
${\cal G}r_j \cap {\cal H}_j$ over $T$.

\begin{definit}$\left.\right.$\vspace{.5pc}

{\rm  \noindent Let $F^{ss}_{\bar L}$ be the functor
$F^{ss}_{\bar L}\hbox{\rm :}\ \hbox{Schemes} \rightarrow
\hbox{Sets}$ which associates to a scheme $T$ the set of
equivalence classes of families $({\cal E}, F({\cal E}))
\rightarrow T\times X$ of $\alpha$-semistable GPBs of rank $r$ and
degree $d$ with $\det {\cal E}_t \cong L$ for all $t\in T$ such
that the section of ${\cal G}r$ defined by $(F_j({\cal E})_j)$
maps into $H_T$.

One similarly defines a full subfunctor $F^s_{\bar L}$ of
$F^{ss}_{\bar L}$ with semistable replaced by stable.}
\end{definit}

\subsection{\it Construction of the moduli space}

Let $S$ denote the set of semistable GPBs $(E, F(E))$, where $E$
is a vector bundle of rank $r$, degree $d$ with fixed determinant
$L$ and $F(E)$ is a subspace of $E_{x_j}\oplus E_{x'_j}$ of
dimension $r$ with fixed weights $(0, \alpha),  0<\alpha <1$. For
$m \gg 0$, all GPBs in $S$ satisfy the condition
\begin{align*}
\hskip -4pc (*) \hskip 3pc H^1(E(m))= 0, H^0(E(m)) \cong {\bf C}^n, H^0(E(m))\!
\rightarrow\!H^0(E(m)\otimes (\oplus_j   {\cal O}_{D_j}))
\end{align*}
is surjective.

Let $Q$ denote the quot scheme of coherent quotients of ${\cal
O}_X^n\otimes {\cal O}_X(-m)$ with fixed Hilbert polynomial
determined by $r, d$. Let $R$ be the nonsingular subvariety of $Q$
corresponding to quotient vector bundles $E$ satisfying condition
$(*)$. Denote by $R_0$ the nonsingular closed subvariety in $R$
corresponding to $E$ with $\Lambda^r (E)= L$. Let ${\cal E}
\rightarrow R\times X$ be the universal quotient bundle. Over $R$,
we have the fibre bundle ${\cal G}r$ with each fibre an
$m$-fold product of $\hbox{Gr}(r,2r)$ as in \S2.2. Over $R_0$, we
have the fibre bundle $\tilde {R_0} := H_{R_0}$ whose fibres are
$m$-fold products of hyperplane sections of $\hbox{Gr}(r,2r)$.
Then $\tilde {R_0}$ is a closed subvariety of ${\cal G}r$. Let
$\tilde {R_0}^{\!\!s}$ (resp. $\tilde {R_0}^{\!\!ss}$) be the open set of
stable (resp. semistable) points in $\tilde {R_0}$.

The GIT quotient $M$ of ${\cal G}r$ by PGL$(n)$ for a suitable
polarization (depending on $\alpha$) is the coarse moduli space
for GPBs \cite{1,2}. It is a normal projective variety. Since
$\tilde {R_0}$ is a PGL$(n)$-invariant closed subscheme
(subvariety) of ${\cal G}r$, the GIT quotient $M_{\bar L}$ of
$\tilde{R_0}^{\!\!ss}$ by PGL$(n)$ is a closed subvariety of $M$ (with
a natural reduced subscheme structure). The GIT quotient
$M^s_{\bar L}= \tilde {R_0}^{\!\!s}\!/\!/\hbox{PGL}(n)$ is an open subscheme
of $M_{\bar L}$. It is easy to see that $M_{\bar L}$
(resp. $M^s_{\bar L}$) is the coarse moduli space for the
functor $F^{ss}_{\bar L}$ (resp. $F^s_{\bar L}$).

\setcounter{theore}{0}
\begin{theor}[\!]
Let $\alpha \in (0,1)$. Then there is a coarse moduli space
$M_{\bar L}$ {\rm (}resp. $M^s_{\bar L}${\rm )} for the
functor $F^{ss}_{\bar L}$ {\rm (}resp. $F^s_{\bar
L}${\rm )}. The moduli space $M_{\bar L}$ is a projective
{\rm (}irreducible{\rm )} variety{\rm ,} containing $M^s_{\bar L}$ as an
open subvariety.
\end{theor}

Let $M'$ (resp. $M'_{\bar L}$) be the open subscheme of $M$
(resp. $M_{\bar L}$) consisting of GPBs $(E,F_j(E))$ such
that the projections $p_j, p'_j$ are isomorphisms for all $j$.
Then $M'_{\bar L}$ corresponds to GPBs in $M$ with fixed
determinant ${\bar L}$ with $F_j(L)=(1:\lambda_j),
\lambda_j\in k^*$ and $M_{\bar L}$ is the closure of
$M'_{\bar L}$.

\section{Application to nodal curves}

\renewcommand\thesubsection{\thesection\arabic{subsection}.}

\subsection{}

Let $Y$ be an irreducible projective nodal curve with nodes $y_j,
j= 1, \dots, m$ and $X$ its desingularization with $D_{\!j}$, the
inverse image of $y_j$. Then there is a correspondence from GPBs
on $X$ of rank $r$, degree $d$ to torsion-free sheaves on $Y$ of
the same rank and degree \cite{1,2}. The correspondence induces a
surjective morphism $f$ from $M$ onto $U$, where $U$ is the moduli
space of torsion-free sheaves of rank $r$, degree $d$ on $Y$. The
restriction of the morphism $f$ to $M'$ is an isomorphism onto the
open subvariety $U'$ of $U$ corresponding to vector bundles on
$Y$. One has $f^{-1} U' = M'$.

For $r=1$, a GPB ${\bar L}$ corresponds to a torsion-free
sheaf ${\cal L}$ on $Y$. Then ${\cal L}$ is a line bundle if and
only if $F_j(L)= (1,\lambda_j), \, \lambda_j\in k^*, \, \forall
j$. Suppose that ${\cal L}$ is a line bundle and let $U'_{\cal L}$
be the closed subset of $U'$ corresponding to vector bundles with
fixed determinant ${\cal L}$. Then $f^{-1} U'_{\cal L}=
M'_{\bar L}$ and the morphism $f$ maps $M'_{\bar L}$
isomorphically onto $U'_{\cal L}$. Hence, if ${\cal L}$ is a line
bundle, then $f(M_{\bar L})$ contains $U'_{\cal L}$ as an
open subset. Since $M_{\bar L}$ is an irreducible, closed
subscheme of  $M$, the image $f(M_{\bar L})$ is a closed,
irreducible subscheme of $U$ and $U'_{\cal L}$, being open, is
dense in it. It follows that $f(M_{\bar L})$ is the closure
of $U'_{\cal L}$.

\begin{rem}
{\rm The projective scheme $M_{\bar L}$ can be regarded as
the compactified moduli space of vector bundles of rank $r$ with
determinant ${\cal L}$ on the nodal curve $Y$.}
\end{rem}

\begin{rem}
{\rm In fact, for any torsion-free sheaf ${\cal L}$ of rank $1$,
the image $f(M_{\bar L})$ is a closed, irreducible subscheme
of $U$ containing the subset of $U$ consisting of torsion-free
sheaves with fixed determinant ${\cal L}$ as an open dense set.}
\end{rem}

\renewcommand\thesubsection{\thesection\arabic{subsection}}

\subsection{\it Relation to Seshadri--Nagaraj conjecture}

For a torsion-free sheaf $F$ of rank $r$ on $Y$, let $N:=
(\Lambda^r F)/\hbox{(torsion)}$, where (torsion) denotes the
maximum subsheaf with proper support. Denote by $I_j$ the ideal
sheaf of the node $y_j$. Define $U_{\cal L}$ as the set
\begin{equation*}
U_{\cal L} = \{ F\in U\hbox{\rm :}\ I_j^r {\cal L} \subset N
\subset {\cal L}, \quad \forall j \}.
\end{equation*}
Seshadri and Nagaraj had defined this set for $Y$ with one node
and conjectured that if ${\cal L}$ is a line bundle, then $U_{\cal
L}$ is the closure of $U'_{\cal L}$ (Conjecture (a), page 136 of
\cite{3}). We prove this conjecture.

Let ${\tilde R}_0^{1-ss}$ denote the subset consisting of
$1$-semistable points, then ${\tilde R}_0^{ss} \subset {\tilde
R}_0^{1-ss}$. Let $P$ be the moduli space of $1$-semistable GPBs
[5]. One has morphisms $f\hbox{\rm :}\ {\tilde R}_0^{ss}
\rightarrow U$ inducing $f\hbox{\rm :}\ M_{\bar L}
\rightarrow U$ and $f_1\hbox{\rm :}\ {\tilde R}_0^{1-ss}
\rightarrow U$ inducing $f_1\hbox{\rm :}\ P \rightarrow U$.

\begin{propo}$\left.\right.$\vspace{.5pc}

\noindent Let ${\bar L}$ be any $G\!P\!B$ of rank $1$ on $X$ and
${\cal L}$ the corresponding torsion-free sheaf of rank $1$ on
$Y$.

\begin{enumerate}
\renewcommand\labelenumi{{\rm (\arabic{enumi})}}
\leftskip .1pc
\item If $(E, F_j(E)) \in {\tilde R}_0^{1-ss}${\rm ,} then $F = f_1((E,F_j(E)))
\in U_{\cal L}$ and hence $f(M_{\bar L})\subset U_{\cal L}$.

\item The morphism ${\tilde R}_0^{1-ss} \rightarrow U$ surjects onto
$U_{\cal L}$.

\item $f(M_{\bar L}) = U_{\cal L}$ for $\alpha$ sufficiently close to $1$.
\end{enumerate}
\end{propo}

\begin{proof}
We may assume that $Y$ has only one node $y$. It is easy to see
(from the proof) that the general case follows exactly on same
lines. Consider a GPB $(E,F(E))$. Let $p_i\hbox{\rm :}\ F(E)
\rightarrow E_{x_i}, \, i= 1, 2$ be the projections and $a_i =
\dim \ker p_i$. Let $E_0= p^*(F)/\hbox{(torsion)}$, then $E_0
\subset E$. Since $F\!\!\mid_{Y - y}\ = p_*E\!\mid_{Y - y}$, one has
$N\!\!\mid_{Y-y}\ = (p_*L)\!\!\mid_{Y-y}\ = {\cal L}\!\mid_{Y-y}$. Hence to
check that $I^r {\cal L} \subset N \subset {\cal L}$, we have only
to check it locally at the node $y$.

Let $(A,m)$ be the local ring at $y$. Its normalization
${\bar A}$ is a semilocal ring with two maximal ideals $m_1,
m_2$. The inclusion
\begin{equation*}
F_y \subset (p_* E)_y
\end{equation*}
may be identified with the inclusion
\begin{equation*}
(r-a_1-a_2)A \oplus a_1 m_1 \oplus a_2 m_2 \subset r {\bar A}
\end{equation*}
(Proposition 4.2 of \cite{2}).

\begin{enumerate}
\renewcommand\labelenumi{(\arabic{enumi})}
\leftskip .1pc
\item We consider the following cases separately.

\begin{case}
{\rm Suppose that $p_1, p_2$ are both isomorphisms. Then $\det (E,
F(E))={\bar L}$ corresponds to a torsion-free sheaf ${\cal
L}$ which is locally free at $y$. In this case, $F$ is locally
free at $y$ with $(\det F)_y = {\cal L}_y$ so that $N={\cal L}
\supset I^r {\cal L}$.}
\end{case}

\begin{case}
{\rm Assume that $p_1$ is an isomorphism, $p_2$ is not an
isomorphism (the opposite case can be dealt with similarly). Then
$\det (E,F(E)) = (L, L_{x_1})= {\bar L}$ corresponds to
${\cal L}$ which is not locally free at $y$. One always has $N
\subset p_*L$ and $N_y \subset {\cal L}_y$ if and only if
$N_y\otimes k(y) \subset F(L)$. Locally, $a_1 =0$ so that $F_y =
(r-a_2)A \oplus a_2 m_2$. Then $N_y = (m_2)^{a_2}$ so that
$N_y\otimes k(y) \subset L_{x_1}= F(L)$. Hence $N \subset {\cal
L}$. Since $m^r \subset m_2^{a_2}$, it follows that $I^r {\cal L}
\subset N \subset {\cal L}$.}
\end{case}

\begin{case}
{\rm If both $p_1$ and $p_2$ are not isomorphisms, one has
$a_1\geq 1, \, a_2\geq 1$. Then locally, $N_y = m_1^{a_1}
m_2^{a_2} \subset m_1 m_2 =m$. It follows that $N_y$ maps to zero
in $p_*L\otimes k(y)$ so that $N \subset I {\cal L} \subset {\cal
L}$. Since $m^r \subset m_1^{a_1} m_2^{a_2}$, one has $I^r {\cal
L}\subset N \subset {\cal L}$.}
\end{case}

\parindent 1pc Note that any $(E,F(E)) \in M_{\bar L}$ with $F(L)=
(1 : \lambda), \lambda \in k^*$ (i.e. ${\cal L}$
locally free at $y$) occurs only in cases (i) or (iii). For
$(E,F(E))\in M_{\bar L}$ with $F(L) = (1:0)$ (or
$F(L)=(0:1)$) only cases (ii) and (iii) occur. Part (1) now
follows. Note that $(E,F_j(E))$ is $1$-semistable if and only if
$F$ is semistable \cite{1,2}.

\item Let $F\in U_{\cal L}$. Since $I^r {\cal L} \subset N \subset
{\cal L}$ it follows that $L\!\!\!\mid_{X-D}\ = p^*N\!\!\!\mid_{X-D}$ where $D=
\sum_j (x_j+x'_j)$. Since $\det E_0 = p^* N$ outside $D$, one has
$L= \det E_0$ outside $D$. It follows that $L = \det E_0\otimes
{\cal O}_X(\sum_j (a_jx_j+a'_jx'_j)), \, a_j+a'_j\leq r$. Let $E$
be given by an extension
\begin{equation*}
0 \rightarrow E_0 \rightarrow E \rightarrow \oplus_j (k(x_j)^{a_j}
\oplus k(x'_j)^{a'_j}) \rightarrow 0.
\end{equation*}
The composite $F \hookrightarrow p_*E_0 \hookrightarrow p_*E$
induces a linear map $F\,\otimes\,k(y_j) \rightarrow p_*(E)\,\otimes\,k(y_j)$.
Let $F_j(E)$ be the image of this linear map. Then
$(E,F_j(E))$ maps to $F$ and it is $1$-semistable as $F$ is
semistable \cite{2}. By construction, $\det E = \det E_0\otimes
{\cal O}_X (\sum_j(a_j x_j + a'_j x'_j)) = L$. It follows that
$(E,F_j(E)) \in {\tilde R}_0^{1-ss}$.

\item Let $P_{\bar L}$ denote the closure of ${\tilde
R}_0^{1-ss}/\hbox{PGL}(n)$ in $P$. For $\alpha$ close to $1$,
there is a surjective birational morphism $\phi\hbox{\rm :}\ M
\rightarrow P$ with $f= f_1\circ \phi$. It maps $M_{\bar L}$
birationally into $P_{\bar L}$. Since both these spaces are
irreducible and of the same dimension, it follows that
$\phi(M_{\bar L})= P_{\bar L}$. Since $f({\tilde
R}_0^{1-ss})$ surjects onto $U_{\cal L}$, it follows that
$f_1(P_{\bar L}) \supset U_{\cal L}$ and hence
$f(M_{\bar L}) \supset U_{\cal L}$. From (1), it follows that
$U_{\cal L}=f(M_{\bar L})= f_1(P_{\bar L})$. \hfill
$\Box$\vspace{-1.5pc}
\end{enumerate}
\end{proof}

\begin{coro}$\left.\right.$\vspace{.5pc}

\noindent If ${\cal L}$ is a line bundle{\rm ,} then $U_{\cal L}$
is the closure of $U'_{\cal L}$.
\end{coro}

\begin{proof}
From Proposition 3.3, one has (as sets) $f(M_{\bar L}) =
U_{\cal L}$. Since  $f(M_{\bar L})$ is the closure of
$U'_{\cal L}$ if ${\cal L}$ is a line bundle, the result
follows.\hfill $\Box$
\end{proof}

\begin{rem}
{\rm The proof of Proposition 3.3(1) easily  generalizes to $(E,
F_j(E))$ replaced by a family $({\cal E}, F_j({\cal
E}))\rightarrow T\times X$.}
\end{rem}

\begin{rem}
{\rm Sun \cite{6} had proved the conjecture (a) of
Seshadri--Nagaraj by considering a smooth curve $X$ degenerating
to an irreducible nodal curve $Y$. However he does not get a
moduli functor or a scheme structure on $U_{\cal L}$ except in
some cases (e.g. one node, rank $2$, degree $1$).}
\end{rem}

\begin{rem}
{\rm Schmitt \cite{4} has constructed a moduli space ${\cal M}$ of
$\alpha$-semistable descending singular SL$(r)$-bundles
$(A,q,\tau)$ where $(A,q)$ is a GPB on $X$ and $\tau\hbox{\rm :}\
$ Sym$^*(A\otimes {\bf C}^r)^{{\rm SL}(r,{\bf C})} \rightarrow
{\cal O}_X$ a nontrivial homomorphism. It is shown that $\det A =
{\cal O}_X$ and for $\alpha \in (0,1)\cap {\bf Q}$, there is a
forgetful morphism $h\hbox{\rm :}\ {\cal M} \rightarrow M$ (\S5.1
of \cite{4}). For $\alpha$ close to $1$, one has a forgetful
morphism ${\cal M} \rightarrow U= U(r,0)$ whose set theoretic
image is $U_{\cal L}$ (Proposition 5.1.1 of \cite{4}). Then, since
$\det A ={\cal O}_X$, it follows that $h({\cal M}) = M_{\bar
L}$ (as\break sets).}
\end{rem}

\section*{Acknowledgements}

I would like to thank I Biswas and A Schmitt for useful comments
and careful reading of the previous version.


\begin{thebibliography}{9}

\bibitem{1} Bhosle Usha N, Generalised parabolic bundles and applications
to torsion-free sheaves on nodal curves. {\it Arkiv for Matematik}
{\bf 30(2)} (1992) 187--215

\bibitem{2} Bhosle Usha N, Generalised parabolic bundles and
applications~II. {\it Proc. Indian Acad. Sci. (Math. Sci.)} {\bf
106(4)} (1996) 403--420

\bibitem{3} Nagaraj D S and Seshadri C S, Degenerations of the moduli
spaces of vector bundles on curves I. {\it Proc. Indian Acad. Sci.
(Math. Sci.)} {\bf 107(2)} (1997) 101--137

\bibitem{4} Schmitt A, Singular principal $G$-bundles on nodal curves.
{\it J. Eur. Math. Soc.} {\bf 7} (2005) 215--251

\bibitem{5} Sun X, Degeneration of moduli spaces and generalized theta functions.
{\it J. Algebraic Geom.} {\bf 9} (2000) 459--527

\bibitem{6} Sun X, Moduli spaces of SL$(r)$-bundles on singular
irreducible curves. {\it Asian J. Math.} {\bf 7(4)} (2003)
609--625, math.AG/0303198
\end{thebibliography}
\end{document}